\documentclass[11pt, reqno]{amsart}
\usepackage{amssymb,amsmath,amscd,amsfonts,amsthm,mathrsfs}
\usepackage[french, english]{babel}
\usepackage{verbatim}
\usepackage{graphicx}             
\input xy
\xyoption{all}
\usepackage[all]{xy}
\usepackage{color}

\usepackage[latin1]{inputenc}

\newcommand{\N}{\mathbb{N}}
\newcommand{\Z}{\mathbb{Z}}

\newcommand{\R}{\mathbb{R}}
\newcommand{\C}{\mathbb{C}}
\newcommand{\Ac}{\mathcal{A}}

\newcommand{\Cc}{\mathcal{C}}
\newcommand{\Ec}{\mathcal{E}}
\newcommand{\Fc}{\mathcal{F}}
\newcommand{\Ic}{\mathcal{I}}
\newcommand{\Jc}{\mathcal{J}}
\newcommand{\Kc}{\mathcal{K}}

\newcommand{\Hc}{\mathcal{H}}
\newcommand{\Lc}{\mathcal{L}}

\newcommand{\Tc}{\mathcal{T}}
\newcommand{\Vc}{\mathcal{V}}

\def\build#1_#2^#3{\mathrel{\mathop{\kern 0pt#1}\limits_{#2}^{#3}}}

\newtheorem{theorem}{Theorem}[section]
\newtheorem{propr}[theorem]{Proposition}
\newtheorem{lem}[theorem]{Lemma}
\newtheorem{rem}[theorem]{Remark}
\newtheorem{exa}[theorem]{Example}
\newtheorem{defi}[theorem]{Definition}
\newtheorem{coro}[theorem]{Corollary}
\newtheorem{thm}[theorem]{Theorem}
\newenvironment{dem}{\noindent {\it Proof: }
    \begin{quotation}\noindent}{\end{quotation}\hfill$\square$}
\newenvironment{demL}{\noindent {\it Proof of Lemma 4.11: }
    \begin{quotation}\noindent}{\end{quotation}\hfill$\square$}
\newenvironment{demP}{\noindent {\it Proof of Proposition 4.10: }
    \begin{quotation}\noindent}{\end{quotation}\hfill$\square$}

\numberwithin{equation}{section}

\title[Spectral Gap of the Laplacian on triangulations]
{Spectral Gap of The Discrete Laplacian \\on Triangulations}

\begin{document}
\title[Spectral Gap of the Laplacian on triangulations]
{Spectral Gap of The Discrete Laplacian \\on triangulations}

\author{Yassin CHEBBI}
\address{Laboratoire de Math\'ematiques Jean Leray, Facult\'{e} des Sciences, CNRS,
Universit\'{e} de Nantes, BP 92208, 44322 Nantes, France}
\address{Universit\'e de Monastir, LR/18ES15, Tunisie}
\email{chebbiyassin88@gmail.com,~yassin.chebbi@univ-nantes.fr}
\date{Version of \today}

\begin{abstract}
Our goal in this paper is to find an estimate for the spectral gap
of the Laplacian on a 2-simplicial complex consisting on a triangulation
of a complete graph. An upper estimate is given by generalizing the Cheeger
constant. The lower estimate is obtained from the first non-zero eigenvalue
of the discrete Laplacian acting on the functions of certain sub-graphs.

\end{abstract}


\maketitle

\tableofcontents

\maketitle

\section{Introduction}
The concept of triangulation was investigated in \cite{Che} as a generalization of graphs.
This structure of a 2-simplicial complex allows to define our discrete Laplacian which
acts on the triplets of functions, 1-forms and 2-forms. This paper deals with questions
on spectral theory of triangulations for Laplacians. There are several recent works giving
lower bounds of the spectrum of the Laplacian on graphs via isoperimetric estimates,
see (\cite{AM}, \cite{D} and \cite{G}). Moreover, Cheeger gave in \cite{Ch} estimates
of the first non zero eigenvalue of the Laplace-Beltrami operator on a compact manifold
in terms of a geometric constant. This inspired a similar theory on graphs for the Laplacian
acting on the functions, see (\cite{Chu}, \cite{CGY}, \cite{CdV} and \cite{KL}).

Considering first only locally finite graphs, we use the Weyl criterion, known from \cite{RS}, to show that all
the spectra of our different Laplacians are in connection except for the value $0.$
Starting from Section 3, we restrict our work to finite triangulations.
Our main interest is about the minimal eigenvalue of the upper Laplacian $\Lc^+_1=\delta^1d^1$
on 1-forms, to be called the spectral gap. Since $Im(d^0)\subset\ker(\Lc_1^+)=\ker(d^1),$
there are $dim(Im(d^0))=|\Vc|-1$ trivial $0$-eigenvalues, if $|\Vc|$ is the number
of vertices and we define the spectral gap to coincide with the $|\Vc|^{th}$-eigenvalue
of $\Lc_1^+.$ More precisely, we discuss lower and upper estimates for the spectral gap.
It is an open question on the Riemannian manifold case. An upper estimate is given by
generalizing Cheeger's approach for triangulations of a complete graph $\Tc_n.$
The Cheeger constant is defined as follows, see \cite{Chee} and \cite{PRT}:
$$h(\Tc_n):=\displaystyle\min_{\displaystyle\Vc=\displaystyle\bigcup_{i=0}^{2}\Ac_i}
\dfrac{n|\Fc(\Ac_0,\Ac_1,\Ac_2)|}{|\Ac_0||\Ac_1||\Ac_2|},
$$
where $\Ac_0,\Ac_1,\Ac_2$ are nonempty sets, making a partition of $\Vc,$
and $\Fc(\Ac_0,\Ac_1,\Ac_2)$ denotes the set of the triangle faces with
one vertex in each $\Ac_i.$

Moreover, we obtain a lower bound in terms of the first non-zero eigenvalue of
the discrete Laplacian $\Lc_0$ defined on the space of functions on the vertices
of certain sub-graphs.

\section{Preliminaries}
\subsection{Notion of graph}
Let $\Vc$ be a countable set of \emph{vertices} and $\Ec$ a subset of $\Vc\times\Vc,$
the set of \emph{oriented edges}. The pair $\Kc=$($\Vc,\Ec$) is called a \emph{graph}.
We assume that $\Ec$ is \emph{symmetric}, \emph{ie.} $(x,y)\in\Ec\Rightarrow(y,x)\in\Ec.$
When two vertices $x$ and $y$ are connected by an edge $e$,
we say they are \emph{neighbors}. We denote $x\sim y$ and $e=(x,y)\in\Ec.$
The set of neighbors of $x\in\Vc$ is denoted by $\Vc(x):=\{y\in\Vc,~y\sim x\}.$
The graph $\Kc$ is said \emph{locally finite} if each vertex belongs to a finite
number of edges. The \emph{degree} or \emph{valence} of a vertex $x\in\Vc$
is the cardinal of the set $\Vc(x),$ denoted by $deg(x).$ If the graph $\Kc$ has
a finite set of vertices, it is called a \emph{finite graph}. An \emph{oriented graph}
$\Kc$ is given by a partition of $\Ec:$
$$\Ec=\Ec^-\cup\Ec^+
$$
$$(x,y)\in\Ec^-\Leftrightarrow(y,x)\in\Ec^+.
$$

In this case for $e=(x,y)\in\Ec^-,$ we define the origin $e^-=x,$
the ending $e^+=y$ and the opposite edge $-e=(y,x)\in\Ec^+.$
A \emph{path} is a finite sequence of edges $\{e_i\}_{0\leq i\leq n}$
such that if $n\geq 2$ then $e_i^+=e_{i+1}^-,$ for all $i\in\{0,...,n-1\}.$
The graph $\Kc$ is said \emph{connected} if any two vertices $x$ and $y$
can be connected by a path with $e_0^-=x,~e_n^+=y.$ A \emph{cycle} is a path
whose origin and end are identical, i.e $e_0^-=e_n^+.$ \emph{An n-cycle}
is a cycle with $n$ vertices. If no cycles appear more than once in a path,
the path is called a \emph{simple path}. All the graphs we shall consider on the sequel will be:
\begin{center}
\textbf{connected, oriented, without loops and locally finite.}
\end{center}
\subsection{Notion of triangulation}
A triangulation generalize the notion of a graph. This structure gives
a general framework for Laplacians defined in terms of the combinatorial
structure of a simplicial complex. We refer to (\cite{Chee},\cite{FP}) for more detail.

\emph{The set of direct permutations} of $(x,y,z)\in\Vc^{3}$ is denoted by
$$\sigma(x,y,z):=\left\{(x,y,z),(y,z,x),(z,x,y)\right\}.
$$

Let $Tr$ be the set of all simple 3-cycles. We denote $\Fc$ the set of
\emph{triangular faces} which is a subset of $Tr$ quotiented by direct permutations as follows:
$$\Fc:=Tr/\simeq
$$
where $\varpi_{1}\simeq\varpi_{2}$ if $\varpi_{1}$ is a direct permutation of $\varpi_{2}.$

\begin{defi}
A \emph{triangulation} $\Tc$ is the triplet $(\Vc,\Ec,\Fc),$
where $\Vc$ is the set of vertices, $\Ec$ is the set of edges
and $\Fc$ is the set of \emph{triangular faces}. This structure
is denoted also by the pair $(\Kc,\Fc),$ where $\Kc=(\Vc,\Ec)$
is a connected locally finite graph. Indeed, one can have simple
3-cycles that are not oriented faces. A triangulation is said complete,
if all the triangles are faces.
\end{defi}
\begin{rem}
In this paper, the triangulations are considered as two-dimensional simplicial complexes
where all faces are triangles.
\end{rem}
\begin{figure}[!ht]
\centering
\begin{minipage}[t]{10cm}
\centering
\includegraphics[width=13cm,height=6cm]{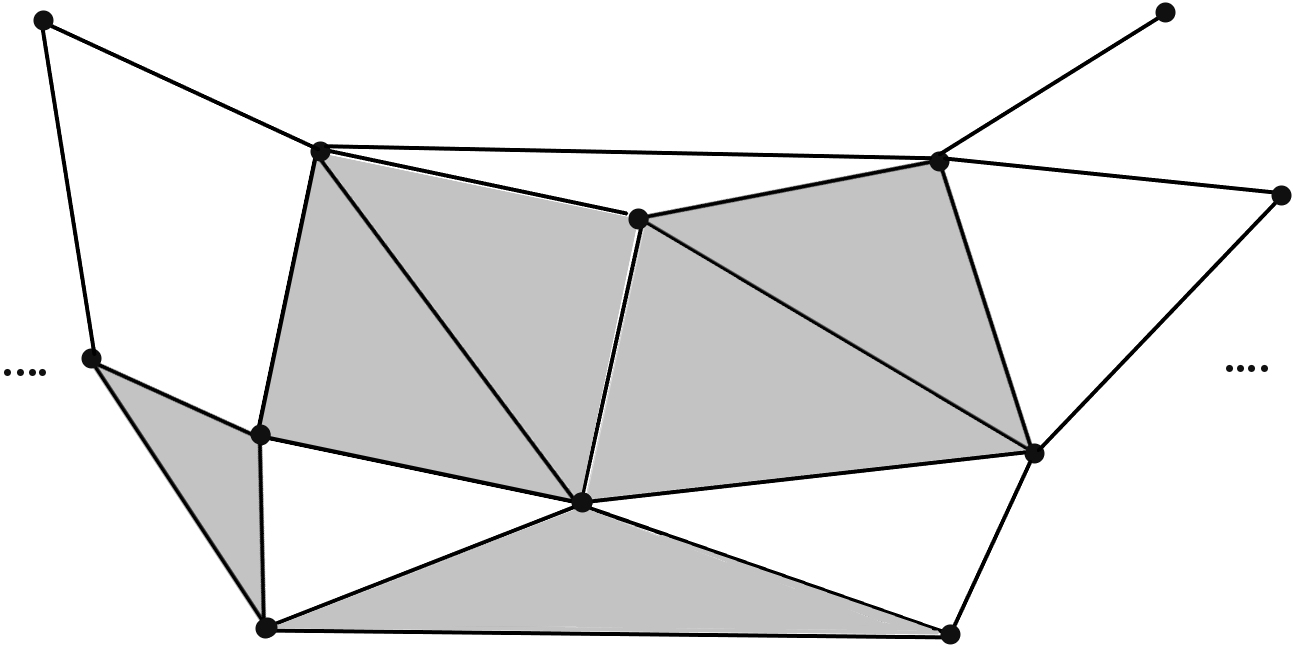}
\caption{Triangulation}
\end{minipage}
\end{figure}

Choosing an orientation of a triangulation consists of defining a partition of $\Fc:$
$$\Fc=\Fc^{-}\cup\Fc^{+},
$$
$$(x,y,z)\in\Fc^{-}\Longleftrightarrow(y,x,z)\in\Fc^{+}.
$$

For each face $\varpi=(x,y,z)\in\Fc,$ we denote $-\varpi$ the opposite face $(y,x,z):$
$$\varpi=(x,y,z)=(y,z,x)=(z,x,y)\in\Fc\Longleftrightarrow-\varpi=(y,x,z)=(x,z,y)=(z,y,x)\in\Fc.
$$

Let $\Tc=(\Vc,\Ec,\Fc)$ be a triangulation. We say that $\Tc$ is a triangulation of \emph{bounded degree},
if the exist a non-negative $C$ such that for all $x\in\Vc$ we have $deg(x)<C.$ For an edge $e\in\Ec$,
we also denote the oriented face $(e^-,e^+,x)$ by $(e,x)$ and denote $\Fc_e$ the set of neighbors of the edge $e:$
$$\Fc_e:=\{x\in\Vc,~(e,x)\in\Fc\}\subseteq\Vc(e^-)\cap\Vc(e^+).
$$

To define weighted triangulations we need weights, let us give
\begin{itemize}
\item $c: \Vc\rightarrow\R^*_+$ the weight on the vertices.
\item $r:\Ec\rightarrow\R^*_+$ the even weight on the oriented edges,
\emph{i.e} $\forall e\in\Ec,~r(-e)=r(e).$
\item $s:\Fc\rightarrow\R^*_+$ the even weight on oriented faces,
\emph{ie.} $\forall \varpi\in\Fc,~s(-\varpi)=s(\varpi).$
\end{itemize}

The weighted triangulation $(\Tc,c,r,s)$ is given by the triangulation $\Tc=(\Kc,\Fc).$
A triangulation $\Tc$ is called \emph{homogeneous}, if the weights of the vertices,
the edges and faces equal $1.$ As our graph $\Kc$ is locally finite, we can define a weight on $\Vc$ by
$$\widetilde{c}(x)=\displaystyle\sum_{y\sim x}r(x,y),~x\in\Vc.
$$

This weight $\widetilde{c}$ on any vertex is well defined. A graph $\Kc$
is called a \emph{normalized graph} if $c(x)=\widetilde{c}(x),$ for all $x\in\Vc.$

\subsection{Functional spaces}
We denote the set of 0-cochains or functions on $\Vc$ by:
$$\Cc(\Vc)=\left\{f:\Vc\rightarrow\C\right\}
$$
and the set of functions of finite support by $\Cc_{c}(\Vc).$ Similarly, we denote the set
of 1-cochains or 1-forms on $\Ec$ by:
$$\Cc(\Ec)=\{\varphi:\Ec\rightarrow\C,~\varphi(-e)=-\varphi(e)\}
$$
and the set of 1-forms of finite support by $\Cc_{c}(\Ec).$ Moreover, we denote the set of
2-cochains or 2-forms on $\Fc$ by:
$$\Cc(\Fc)=\{\phi:\Fc\rightarrow\C,~\phi(-\varpi)=-\phi(\varpi)\}
$$
and the set of 2-forms of finite support by $\Cc_{c}(\Fc).$

We consider on the weighted triangulation $(\Tc,c,r,s)$ the following Hilbert spaces, let us give
\begin{itemize}
\item The  Hilbert space of functions:
$$l^2(\Vc):=\left\{f\in\Cc(\Vc);~\displaystyle\sum_{x\in\Vc}c(x)|f(x)|^{2}<\infty\right\}
$$
with the inner product
$$\langle f,g \rangle_{l^2(\Vc)}:=\displaystyle\sum_{x\in\Vc}c(x)f(x)\overline{g(x)}.
$$
\item The  Hilbert space of 1-forms:
$$l^{2}(\Ec):=\left\{\varphi\in\Cc(\Ec);
~\displaystyle\sum_{e\in\Ec}r(e)|\varphi(e)|^{2}<\infty\right\}
$$
with the inner product
$$\langle \varphi,\psi \rangle_{l^{2}(\Ec)}
=\dfrac{1}{2}\displaystyle\sum_{e\in\Ec}r(e)\varphi(e)\overline{\psi(e)}
$$
\item The Hilbert space of 2-forms:
$$l^{2}(\Fc):=\left\{\phi\in\Cc(\Fc);
~\displaystyle\sum_{(x,y,z)\in\Fc}s(x,y,z)|\phi(x,y,z)|^{2}<\infty\right\}
$$
with the inner product
$$\langle \phi,\theta \rangle_{l^{2}(\Fc)}
=\dfrac{1}{2}\displaystyle\sum_{(x,y,z)\in\Fc}s(x,y,z)\phi(x,y,z)\overline{\theta(x,y,z)}
$$
\end{itemize}

The direct sum of the spaces $l^{2}(\Vc)$, $l^{2}(\Ec)$ and $l^{2}(\Fc)$
can be considered as a Hilbert space denoted by $\Hc$, that is
$$\Hc:=l^{2}(\Vc)\oplus l^{2}(\Ec)\oplus l^{2}(\Fc),
$$
endowed with the inner product
$$\langle (f_{1},\varphi_{1},\phi_{1}),(f_{2},\varphi_{2},\phi_{2}) \rangle_{\Hc}
=\langle f_{1},f_{2} \rangle_{l^2(\Vc)}
+\langle \varphi_{1},\varphi_{2} \rangle_{l^{2}(\Ec)}
+\langle \phi_{1},\phi_{2} \rangle_{l^{2}(\Fc)}.
$$
\subsection{Operators}
In this section, we recall the concept of difference and exterior derivative
operators introduced on weighted triangulations. We refer to \cite{AT} and
\cite{Chee} for more details. This permits to define the discrete Laplacians
acting on functions, 1-forms and 2-forms.

\subsubsection{The difference operator}
It is the operator $d^{0}:\Cc_{c}(\Vc)\longrightarrow \Cc_{c}(\Ec),$ given by
  $$d^{0}(f)(e)=f(e^{+})-f(e^{-}).
  $$

\subsubsection{The co-boundary operator}
It is the formal adjoint of $d^0,$  denoted by $\delta^0.$
$$\delta^0:\Cc_{c}(\Ec)\longrightarrow \Cc_{c}(\Vc),
$$
acts as
$$\delta^{0}(\varphi)(x)=\dfrac{1}{c(x)}\displaystyle\sum_{e,e^{+}=x}r(e)\varphi(e),
$$
and satisfies
\begin{align*}
  \langle d^{0}f,\varphi\rangle_{l^{2}(\Ec)}
  & =\langle f,\delta^{0}\varphi\rangle_{l^{2}(\Vc)},~\forall (f,\varphi)\in\Cc_{c}(\Vc)\times\Cc_{c}(\Ec).
\end{align*}

\subsubsection{The exterior derivative}
It is the operator $d^{1}:\Cc_{c}(\Ec)\longrightarrow \Cc_{c}(\Fc),$ given by
$$d^{1}(\varphi)(x,y,z)=\varphi(x,y)+\varphi(y,z)+\varphi(z,x).
$$

\subsubsection{The co-exterior derivative}
It is the formal adjoint of $d^1,$  denoted by $\delta^{1}.$
$$\delta^1:\Cc_{c}(\Fc)\longrightarrow \Cc_{c}(\Ec),
$$
acts as
$$\delta^{1}(\phi)(e)=\dfrac{1}{r(e)}\sum_{x\in \Fc_{e}}s(e,x)\phi(e,x),
$$
and satisfies
\begin{align*}
  \langle d^{1}\varphi,\phi\rangle_{l^{2}(\Fc)}
  & =\langle \varphi,\delta^{1}\phi\rangle_{l^{2}(\Ec)},~\forall (\varphi,\phi)\in\Cc_{c}(\Ec)\times\Cc_{c}(\Fc).
\end{align*}

\subsubsection{The discrete Laplacian} In this section, we will always consider
a weighted triangulation $(\Tc,c,r,s).$ The discrete Laplacian acting on functions is given by
$$\Lc_{0}(f)(x):=\delta^0d^{0}(f)(x)=\dfrac{1}{c(x)}\displaystyle\sum_{e,e^{+}=x}r(e)d^{0}(f)(e).
$$
for all $f\in\Cc_{c}(\Vc)$ and $x\in\Vc.$ The discrete Laplacian acting on 1-forms is given by
\begin{equation*}
\begin{split}
\Lc_{1}(\varphi)(x,y)
& :=(d^{0}\delta^0+\delta^{1}d^{1})(\varphi)(x,y)\\
&=\dfrac{1}{c(y)}\displaystyle\sum_{e,e^{+}=y}r(e)\varphi(e)
-\dfrac{1}{c(x)}\displaystyle\sum_{e,e^{+}=x}r(e)\varphi(e)\\
&+\dfrac{1}{r(x,y)}\displaystyle\sum_{z\in\Fc_{(x,y)}}s(x,y,z)d^{1}(\varphi)(x,y,z).
\end{split}
\end{equation*}
for all $\varphi\in\Cc_{c}(\Ec)$ and $(x,y)\in\Ec.$ The discrete Laplacian acting on 2-forms is given by
\begin{equation*}
\begin{split}
\Lc_{2}(\phi)(x,y,z)
&:=d^{1}\delta^{1}(\phi)(x,y,z)\\
&=\dfrac{1}{r(x,y)}\displaystyle\sum_{u\in\Fc_{(x,y)}}s(x,y,u)\phi(x,y,u)\\
&+\dfrac{1}{r(y,z)}\displaystyle\sum_{u\in\Fc_{(y,z)}}s(y,z,u)\phi(y,z,u)\\
&+\dfrac{1}{r(z,x)}\displaystyle\sum_{u\in\Fc_{(z,x)}}s(z,x,u)\phi(z,x,u).
\end{split}
\end{equation*}
for all $\phi\in\Cc_{c}(\Fc)$ and $(x,y,z)\in\Fc.$ To define \emph{the Gau\ss-Bonnet} operator,
let us begin by defining the operator
$$d:\Cc_{c}(\Vc)\oplus\Cc_{c}(\Ec)\oplus\Cc_{c}(\Fc)\circlearrowleft
$$
by
$$d(f,\varphi,\phi):=(0,d^{0}f,d^{1}\varphi)
$$
and $\delta$ the formal adjoint of $d$ is given by
$$\delta(f,\varphi,\phi)=(\delta^0\varphi,\delta^{1}\phi,0),~
\forall(f,\varphi,\phi)\in\Cc_{c}(\Vc)\oplus\Cc_{c}(\Ec)\oplus\Cc_{c}(\Fc).
$$

\emph{The Gau\ss-Bonnet} operator is defined on $\Cc_{c}(\Vc)\oplus\Cc_{c}(\Ec)\oplus\Cc_{c}(\Fc)$ into itself by:
\begin{center}
$T:=d+\delta\cong\begin{pmatrix}
0&\delta^0&0 \\
d^{0}&0&\delta^{1}\\
0&d^{1}&0
\end{pmatrix}$
\end{center}

This operator is of Dirac type and is motived by the Hodge Laplacian:
$$\Lc:=T^{2}=\left(d+\delta\right)^{2}=\Lc_{0}\oplus\Lc_{1}\oplus\Lc_{2},
$$

\begin{rem}
The operator $\Lc_1$ is called the \emph{full Laplacian} and defined as $\Lc_1=\Lc_1^-+\Lc_1^+,$
where $\Lc_1^-=d^{0}\delta^0$ \emph{(}resp. $\Lc_{1}^{+}=\delta^{1}d^{1}$\emph{)}
is called \emph{the lower Laplacian} \emph{(}resp. \emph{the upper Laplacian}\emph{)}.
\end{rem}

\section{The spectrum of the Laplacians}
In this section, we will prove the relation between the spectrum of $\Lc^+_1$ and that of $\Lc_2.$
The following results due to \emph{the Weyl's criterion} known from \cite{RS} to characterize
the spectrum of our operators.

\underline{\emph{\textbf{Weyl's criterion:}}} Let $\Hc$ be a separable Hilbert space,
and let $\Lc$ be a bounded self-adjoint operator on $\Hc.$ Then $\lambda$ is in
the spectrum of $\Lc,$ if and only if, there exists a sequence $\left(f_{n}\right)_{n\in\N}$
so that $\|f_{n}\|=1$ and $\displaystyle\lim_{n\longrightarrow\infty}\|\left(\Lc-\lambda\right)f_{n}\|=0.$

We denote $\sigma(\Lc)$ the spectrum of $\Lc.$ We refer here to \cite{Ay}, which proves that $\sigma(\Lc_0)\backslash\{0\}=\sigma(\Lc^-_1)\backslash\{0\}$ in a normalized graph.

\begin{propr}\label{pro1}
Let $\Tc$ be a homogeneous triangulation of bounded degree. Then, we have

$$\sigma(\Lc^+_1)\backslash\{0\}=\sigma(\Lc_2)\backslash\{0\}.
$$
\end{propr}
\begin{dem}
At first, we notice that the hypothesis assures that all operators are bounded.
Through the Weyl's criterion, set $\lambda\in\R^*$ in the spectrum of $\Lc^+_1,$
then there is a sequence $(\psi_n)_n$ in $l^2(\Ec)$ such that:
$$\|\psi_n\|_{l^2(\Ec)}=1\mbox{ and }\displaystyle\lim_{n\rightarrow\infty}\|(\Lc^+_1-\lambda)\psi_n\|_{l^2(\Ec)}=0.
$$

So, we should find a sequence $(\phi_n)_n$ in $l^2(\Fc)$ such that
$$\mbox{ for each }n,~\|\phi_n\|_{l^2(\Fc)}=1\mbox{ and }\displaystyle\lim_{n\rightarrow\infty}\|(\Lc_2-\lambda)\phi_n\|_{l^2(\Fc)}=0.
$$

Set
$$\phi_n:=\dfrac{d^1\psi_n}{\|d^1\psi_n\|_{l^2(\Fc)}},~~n\in\N.
$$

First, let us show that $\|d^1\psi_n\|_{l^2(\Fc)}\neq0.$ We have that
\begin{equation*}
\begin{split}
\|d^1\psi_n\|^2_{l^2(\Fc)}&=\langle \Lc^+_1\psi_n,\psi_n\rangle_{l^2(\Ec)}\\
&=\langle (\Lc^+_1-\lambda)\psi_n,\psi_n\rangle_{l^2(\Ec)}+\langle \lambda\psi_n,\psi_n\rangle_{l^2(\Ec)}\\
&=\langle (\Lc^+_1-\lambda)\psi_n,\psi_n\rangle_{l^2(\Ec)}+\lambda.
\end{split}
\end{equation*}

Thus the term $\|d^1\psi_n\|^2_{l^2(\Fc)}$ tends to $\lambda$ as $n\rightarrow\infty.$
Therefore, there exist $\varepsilon>0$ and $n_0\in\N$ such that for all $n\geq n_0$
we have $\|d^1\psi_n\|_{l^2(\Fc)}>\varepsilon.$ It remains to show that
$\displaystyle\lim_{n\rightarrow\infty}\|(\Lc_2-\lambda)\phi_n\|_{l^2(\Fc)}=0.$ In fact, we have:
\begin{equation*}
\begin{split}
\|(\Lc_2-\lambda)\phi_n\|_{l^2(\Fc)}&=\|(\Lc_2-\lambda)\dfrac{d^1\psi_n}{\|d^1\psi_n\|_{l^2(\Fc)}}\|_{l^2(\Fc)}\\
&=\dfrac{\|(\Lc_2-\lambda)d^1\psi_n\|_{l^2(\Fc)}}{\|d^1\psi_n\|_{l^2(\Fc)}}\\
&=\dfrac{\|\Lc_2d^1-\lambda d^1\psi_n\|_{l^2(\Fc)}}{\|d^1\psi_n\|_{l^2(\Fc)}}\\
&=\dfrac{\|d^1(\Lc^+_1-\lambda)\psi_n\|_{l^2(\Fc)}}{\|d^1\psi_n\|_{l^2(\Fc)}}.
\end{split}
\end{equation*}

Since the operator $d^1$ is bounded. Hence, there exists a constant $c>0$ such that $\|d^1\|\leq c$ and
we have that:
$$\|(\Lc_2-\lambda)\phi_n\|_{l^2(\Fc)}\leq\dfrac{c}{\varepsilon}\|(\Lc^+_1-\lambda)\psi_n\|_{l^2(\Ec)}.
$$

Using the same method as the first step, set the non-zero constant $\lambda\in\sigma(\Lc_2).$
Then there exists a sequence $(\phi_n)$ in $l^2(\Fc)$ such that $\|\phi_n\|_{l^2(\Fc)}=1.$
We consider a sequence $(\psi_n)$ in $l^2(\Ec),$ define as follows:
$$\psi_n:=\dfrac{\delta^1\phi_n}{\|\delta^1\phi_n\|_{l^2(\Fc)}}.
$$
\end{dem}
\begin{propr}\label{pr}
Let $\Tc$ be a finite weighted triangulation. Then, we have
$$\min\sigma(\Lc_1)\leq\displaystyle\min_{e\in\Ec}
\left[\left(\dfrac{1}{c(e^-)}+\dfrac{1}{c(e^+)}\right)r(e)+deg_{\Ec}(e)\right].
$$
\end{propr}
\begin{dem}

Let $e\in\Ec,$ we consider $\chi^{e}=1_{\{e\}}-1_{\{-e\}}.$ Thus, we have
$$\left\|\dfrac{\chi^{e}}{\sqrt{r(e)}}\right\|_{l^2(\Ec)}=1.
$$

Moreover, we have
\begin{equation*}
\begin{split}
\left\langle\dfrac{\chi^{e}}{\sqrt{r(e)}},\Lc_1\left(\dfrac{\chi^{e}}{\sqrt{r(e)}}\right)\right\rangle_{l^2(\Ec)}
&=\left\langle\dfrac{\chi^{e}}{\sqrt{r(e)}},\Lc^-_1\left(\dfrac{\chi^{e}}{\sqrt{r(e)}}\right)\right\rangle_{l^2(\Ec)}
+\left\langle\dfrac{\chi^{e}}{\sqrt{r(e)}},\Lc^+_1\left(\dfrac{\chi^{e}}{\sqrt{r(e)}}\right)\right\rangle_{l^2(\Ec)}\\
&=\left(\dfrac{1}{c(e^-)}+\dfrac{1}{c(e^+)}\right)r(e)+deg_{\Ec}(e).
\end{split}
\end{equation*}

By Rayleigh principle, we obtain that
\begin{equation*}
\begin{split}
\min\sigma(\Lc_1)&\leq
\displaystyle\min_{e\in\Ec}\left[\dfrac{\left\langle\dfrac{\chi^{e}}{\sqrt{r(e)}},\Lc_1
\left(\dfrac{\chi^{e}}{\sqrt{r(e)}}\right)\right\rangle_{l^2(\Ec)}}{\left\|\dfrac{\chi^{e}}{\sqrt{r(e)}}
\right\|_{l^2(\Ec)}}\right]\\
&=\displaystyle\min_{e\in\Ec}\left[\left(\dfrac{1}{c(e^-)}+\dfrac{1}{c(e^+)}\right)r(e)
+deg_{\Ec}(e)\right].
\end{split}
\end{equation*}
\end{dem}
\begin{propr}
Let $\Tc$ be a homogeneous finite triangulation. Then, we have
$$\min\sigma(\Lc_2)\leq8+2\displaystyle\min_{e\in\Ec}|\Fc_e|.
$$
\end{propr}
\begin{dem}
Let $e_{0}\in\Ec,$ we consider $\chi^{e_{0}}=1_{\{e_{0}\}}-1_{\{-e_{0}\}}.$
Then, we have
\begin{equation*}
\begin{split}
\left\langle d^1\chi^{e_{0}},\Lc_2d^1\chi^{e_{0}}\right\rangle_{l^2(\Fc)}
&=\left\|\Lc^+_1\chi^{e_{0}}\right\|^2_{l^2(\Ec)}
=\dfrac{1}{2}\displaystyle\sum_{e\in\Ec}\left(\displaystyle\sum_{x\in\Fc_e}d^1\chi^{e_{0}}(e,x)\right)^2
\end{split}
\end{equation*}
Using $(a+b)^2\leq 2(a^2+b^2),$ for all $a,b\in\R.$ We have
\begin{equation*}
\begin{split}
\left\langle d^1\chi^{e_{0}},\Lc_2d^1\chi^{e_{0}}\right\rangle_{l^2(\Fc)}
&\leq\displaystyle\sum_{e\in\Ec}\left(\displaystyle\sum_{x\in\Fc_e}\chi^{e_{0}}(e)\right)^2
+\displaystyle\sum_{e\in\Ec}\left(\displaystyle\sum_{x\in\Fc_e}(\chi^{e_{0}}(e^+,x)+\chi^{e_{0}}(x,e^-))\right)^2\\
&\leq\displaystyle\sum_{e\in\Ec}|\Fc_e|^2\left(\chi^{e_{0}}(e)\right)^2
+4\displaystyle\sum_{e\in\Ec}\left(\displaystyle\sum_{x\in\Fc_e}\chi^{e_{0}}(e^+,x)\right)^2\\
&=2|\Fc_{e_0}|^2+4\displaystyle\sum_{e\in\Ec}\left(\displaystyle\sum_{x\in\Fc_e}
\chi^{e_{0}}(e^+,x)\right)^2\\
&=2|\Fc_{e_0}|^2+4\displaystyle\sum_{x\in\Fc_{e_{0}}}\left(\displaystyle\sum_{e\in\Ec}
\chi^{e_{0}}(e^+,x)\right)^2\\
&=2|\Fc_{e_0}|^2+8|\Fc_{e_0}|.
\end{split}
\end{equation*}

On other hand, we have
$$\|d^1\chi^{e_{0}}\|^2_{l^2(\Fc)}
=\dfrac{1}{2}\displaystyle\sum_{(x,y,z)\in\Fc}|d^1\chi^{e_{0}}(x,y,z)|^2
=\dfrac{1}{2}\left(|\Fc_{e_0}|+|\Fc_{-e_0}|\right)
=|\Fc_{e_0}|.
$$

From Rayleigh principle, we obtain:
$$\min\sigma(\Lc_2)
=\inf_{\phi\in\Cc_{c}(\Fc)\backslash\{0\}}\dfrac{\langle\Lc_2\phi,\phi\rangle_{l^{2}(\Fc)}}
{\|\phi\|_{l^{2}(\Fc)}^{2}}\leq8+2\displaystyle\min_{e\in\Ec}|\Fc_e|.
$$
\end{dem}

\section{Spectral gap of finite triangulation}
\subsection{The spectral gap}
We describe \emph{the discrete Hodge theory} due to Eckmann \cite{Eck}.
This is a discrete analogue of Hodge theory in Riemannian geometry.
Furthermore, it applies to any finite simplicial complex, and not only to manifolds.
Let us begin with
\begin{lem}
We have that
$$\ker(\Lc_1^-)=\ker(\delta^0),~\ker(\Lc_1^+)=\ker(d^1),~\ker(\Lc_2)=\ker(\delta^1)
\mbox{ and }\ker(\Lc_1)=\ker(\delta^0)\cap\ker(d^1).
$$
\end{lem}
\begin{dem}
It is clear that $\ker(\delta^0)\subseteq\ker(\Lc_1^-).$ On the other hand,
if $\Lc_1^-\varphi=d^0\delta^0\varphi=0,$ for all $\varphi\in l^2(\Ec),$ we have that
$$0=\langle \varphi,\Lc_1^-\varphi \rangle_{l^2(\Ec)}=\|\delta^0\varphi\|^2_{l^2(\Vc)}.
$$
Then, $\delta^0\varphi=0,~\forall\varphi\in l^2(\Ec).$ As the previous reasoning, we prove
that $\ker(\Lc_1^+)=\ker(d^1)$ and $\ker(\Lc_2)=\ker(\delta^1).$ After, if $\varphi\in\ker(\Lc_1)$ then
$$0=\langle \varphi,\Lc_1\varphi \rangle_{l^2(\Ec)}=\langle \delta^0\varphi,\delta^0\varphi \rangle_{l^2(\Vc)}
+\langle d^1\varphi,d^1\varphi \rangle_{l^2(\Fc)},
$$
which shows that $\ker(\Lc_1)=\ker(\delta^0)\cap\ker(d^1).$
\end{dem}

Since $\ker(\Lc_1)=\ker(\delta^0)\cap\ker(d^1),$ we have the \emph{discrete Hodge decomposition}
$$l^2(\Ec)=\ker(\Lc_1)\oplus Im(d^0)\oplus Im(\delta^1).
$$

In particular, it follows that the space of harmonic forms can be identified with the homology of $\Tc:$
$$\ker(d^1)/Im(d^0)=Im(\delta^1)^{\perp}/Im(d^0)=\left(Im(d^0)\oplus\ker(\Lc_1)\right)/Im(d^0)\cong\ker(\Lc_1).
$$

The same holds for the homology of $\Tc,$ giving
$$\ker(d^1)/Im(d^0)\cong\ker(\Lc_1)\cong\ker(\delta^0)/Im(\delta^1).
$$

For a triangulation, the space $Im(d^0)$ is always in the kernel of the upper Laplacian,
and considered to be its \emph{trivial zeros}. There can be more zeros in the spectrum, since $Im(d^0)\subset\ker(\Lc_1^+)=\ker(d^1).$ As $l^2(\Ec)=\ker(\delta^0)\oplus Im(d^0),$ this leads
to the following definition:
\begin{defi}\emph{(The spectral gap)}\label{gap}
The spectral gap of a finite triangulation $\Tc,$ denoted $\lambda_{\Tc},$ is the minimal eigenvalue
of the upper Laplacian on 1-forms:
$$\lambda_{\Tc}:=\min\sigma\left(\Lc^+_{1\scriptscriptstyle{\vert\ker(\delta^0)}}\right)
=\min\sigma\left(\Lc_{1\scriptscriptstyle{\vert \ker(\delta^0)}}\right),
$$
the equality follows from $\Lc_{1\scriptscriptstyle{\vert \ker(\delta^0)}}
\equiv \Lc^+_{1\scriptscriptstyle{\vert \ker(\delta^0)}}.$
\end{defi}
\begin{rem}
Generally, the spectral gap is the difference between $0$
and the first non-zero eigenvalue, see \emph{(}\cite{Chu}, \cite{CGY} and \cite{CdV}\emph{)}.
This definition coincides with the definition \ref{gap} in the case
where the eigenvalue $\lambda_{\Tc}$ is not zero \emph{(} i.e there is no harmonic 1-form\emph{)}.
\end{rem}
The following proposition gives a characterization of the spectral gap.
\begin{propr}
The spectrum of the Laplacian $\Lc_1^+$ is composed of $|\Ec|$ reals eigenvalues
$\{\lambda_i\}_{0\leq i\leq|\Ec|-1}$ ranked in ascending order:
$$\lambda_0\leq\lambda_1\dots\leq...\leq\lambda_{|\Ec|-1}.
$$
Moerever, we have that
$$\lambda_{\Tc}=\lambda_{|\Vc|-1}.
$$
\end{propr}
\begin{dem}
Since the dimension of the space $Im(d^0)$ is the number of the trivial zeros of the upper Laplacian.
Using the rank theorem for the operator $d^0:l^2(\Vc)\rightarrow l^2(\Ec),$ we get that
$$\lambda_{\Tc}:=\lambda_{dim(Im(d^0))}=\lambda_{|\Vc|-1}.
$$
\end{dem}
\begin{propr}\label{haifa3}
Let $\Tc$ be a finite triangulation. Then, we have
   $$|\Fc|<|\Ec|-|\Vc|+1\Rightarrow\lambda_{\Tc}=0.
   $$
\end{propr}
\begin{dem}
Applying the Rank theorem for $d^0:l^2(\Vc)\rightarrow l^2(\Ec),$ we obtain
\begin{equation*}
  \begin{split}
     dim(\ker(\delta^0)) & =|\Ec|-dim(Im(d^0))\\
       & = |\Ec|-\left(|\Vc|-1\right)=|\Ec|-|\Vc|+1
  \end{split}
\end{equation*}
And since $dim(Im(\delta^1))\leq|\Fc|,$ then
$$|\Fc|<|\Ec|-|\Vc|+1
\Rightarrow Im(\delta^1)\varsubsetneq\ker(\delta^0)
\Leftrightarrow\lambda_{\Tc}=0.
$$
\end{dem}
\begin{rem}
In the case where $Im(\delta^1)=\ker(\delta^0),$
we have $|\Fc|\geq|\Ec|-|\Vc|+1$ and $\lambda_{\Tc}\neq0.$
\end{rem}
\begin{exa}
Consider the triangulation $\Tc$ such that all the triangles are faces
as in Figure $2.$ Using the Rank theorem, we have:
$$\left\{\begin{array}{ll}
       dim(Im(d^0))=5\hbox{ et }dim(\Cc(\Ec))=12& \Rightarrow dim(\ker(\delta^0))=7\\
       dim(\Cc(\Fc))=8\hbox{ et }dim(\ker(\delta^1))=1& \Rightarrow dim(Im(\delta^1))=7.

                   \end{array}
            \right.
   $$
Then, $\lambda_{\Tc}\neq0.$
\end{exa}
\begin{figure}[!ht]
\centering
\begin{minipage}[t]{10cm}
\centering
\includegraphics*[width=12cm,height=9cm]{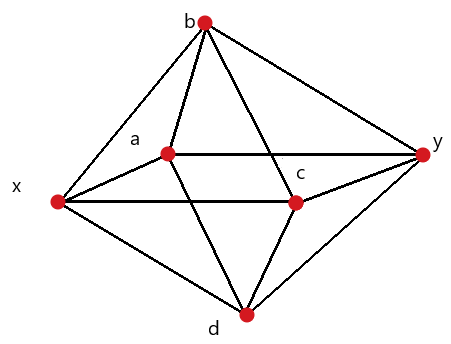}
\caption{Example of a triangulation with non-zero spectral gap }
\end{minipage}
\end{figure}

\subsection{Triangulation of a complete graph}
The aim of this part is to provide effective estimates for the spectral gap of a
triangulation of a complete graph. Let us begin by
\begin{defi}
A triangulation of a complete graph $\Tc_n$ is the pair $(\Kc_n,\Fc),$ where $\Kc_n$
is the complete graph of $n$ vertices and $\Fc$ is the set of triangular faces.
\end{defi}
First, we give a result assuring that our spectral gap is zero in a triangulation
of a complete graph. This result is a consequence of Proposition \ref{haifa3}.
\begin{coro}
Let $\Tc_n$ be a homogeneous triangulation of a complete graph. Then, we have
$$|\Fc|<\dfrac{(n-1)(n-2)}{2}\Rightarrow\lambda_{\Tc}=0.
$$
\end{coro}
\subsubsection{Upper estimates of the spectral gap}
We would like to find a concrete upper estimate of the spectral gap in relation with
the number of oriented faces in a simple triangulation of a complete graph.
\begin{propr}\label{pro0}
Let $\Tc_n$ be a homogeneous triangulation of a complete graph. Then, we have
\begin{enumerate}
\item[\emph{i})]$n\in\sigma(\Lc_1).$
\item[\emph{ii})]If $\Tc_n$ is a simple complete triangulation then $\sigma(\Lc_1)=\{n\}.$
\end{enumerate}
\end{propr}
\begin{lem}\label{lem0}
Let $\Tc$ be a finite triangulation. Then
$$\sigma(\Lc_1)=\sigma(\Lc^-_{1\scriptscriptstyle{\vert Im(d^0)}})
\cup\sigma(\Lc^+_{1\scriptscriptstyle{\vert \ker(\delta^0)}}).
$$
\end{lem}
\begin{demL}
It is clear that $\sigma(\Lc^-_{1\scriptscriptstyle{\vert Im(d^0)}})
\cup\sigma(\Lc^+_{1\scriptscriptstyle{\vert \ker(\delta^0)}})\subseteq\sigma(\Lc_1).$
Let $\lambda\in\sigma(\Lc_1),$ there is $\psi_{\lambda}\in l^2(\Ec)$ such that
$\Lc_1\psi_{\lambda}=\lambda\psi_{\lambda}.$ Since $l^2(\Ec)=\ker(\delta^0)\oplus Im(d^0),$
then  $\psi_{\lambda}=\psi^1_{\lambda}+\psi^2_{\lambda}$ where
$(\psi^1_{\lambda},\psi^2_{\lambda})\in\ker(\delta^0)\times Im(d^0).$ Hence
$$\lambda\in\sigma(\Lc_{1\scriptscriptstyle{\vert Im(d^0)}})
\cup\sigma(\Lc_{1\scriptscriptstyle{\vert \ker(\delta^0)}})
=\sigma(\Lc^-_{1\scriptscriptstyle{\vert Im(d^0)}})
\cup\sigma(\Lc^+_{1\scriptscriptstyle{\vert \ker(\delta^0)}}).
$$
\end{demL}

\begin{demP}
\begin{enumerate}
\item[\emph{i})]Let $\psi\in Im(d^0),$ then it exists $f\in\l^2(\Vc)$ such that $\psi=d^0f.$ We calculate
\begin{equation*}
\begin{split}
\Lc^-_1\psi(x,y)&=\displaystyle\sum_{e,e^+=y}\psi(e)-\displaystyle\sum_{e,e^+=x}\psi(e)\\
&=\displaystyle\sum_{e,e^+=y}\left(f(y)-f(e^-)\right)-\displaystyle\sum_{e,e^+=x}\left(f(x)-f(e^-)\right)\\
&=(n-1)d^0f(x,y)+\displaystyle\sum_{z\sim x}f(z)-\displaystyle\sum_{z\sim y}f(z)\\
&=n\psi(x,y).
\end{split}
\end{equation*}

Hence $\Lc^-_{1\scriptscriptstyle{\vert Im(d^0)}}=nI.$ By Lemma \ref{lem0}, we have that $n\in\sigma(\Lc_1).$
\item[\emph{ii})]Let $\psi\in \ker(\delta^0),$ if $\Tc_n$ is complete, then we have
\begin{equation*}
\begin{split}
\Lc^+_1\psi(x,y)&=\displaystyle\sum_{z\in\Fc_{(x,y)}}\left(\psi(x,y)+\psi(y,z)+\psi(z,x)\right)\\
&=(n-2)\psi(x,y)+\displaystyle\sum_{z\in\Fc_{(x,y)}}\psi(y,z)+\displaystyle\sum_{z\in\Fc_{(x,y)}}\psi(z,x)\\
&=n\psi(x,y)+\delta^0\psi(x)-\delta^0\psi(y)=n\psi(x,y).
\end{split}
\end{equation*}

By Lemma \ref{lem0}, we get that $\sigma(\Lc_1)=\{n\}.$
\end{enumerate}
\end{demP}
\begin{rem}
The value $0$ is not always in the spectrum of the full Laplacian $\Lc_1.$
It depends on the value of $\lambda_{\Tc_n},$ due to the fact that $\min\sigma(\Lc_1)=\min\{\lambda_{\Tc_n},n\}.$
In particular, if $\Tc_n$ is a complete triangulation we have that $\lambda_{\Tc_n}=n.$
\end{rem}
\begin{coro}
Let $\Tc_n$ be a homogeneous triangulation of a complete graph. Then, we have
$$\lambda_{\Tc_n}\leq 2+\displaystyle\min_{e\in\Ec}|\Fc_e|.
$$

Moreover, $\Tc_n$ is the complete triangulation if and only if $\lambda_{\Tc_n}=n.$
\end{coro}
\begin{dem}
From Proposition \ref{pro0} and Lemma \ref{lem0}, we deduce that $\min\sigma(\Lc_1)=\min\{\lambda_{\Tc_n},n\}.$ Moreover, The Proposition \ref{pr}
give an upper estimate of the lower spectrum of $\Lc_1:$
\begin{equation}\label{equa}
\min\sigma(\Lc_1)\leq 2+\displaystyle\min_{e\in\Ec}|\Fc_e|.
\end{equation}

If $\min\sigma(\Lc_1)<n$ then $\lambda_{\Tc_n}=\min\sigma(\Lc_1).$ Then,
the inegality (\ref{equa}) give an upper estimate of the lower spectrum
$\lambda_{\Tc_n}$ acting as:
$$\lambda_{\Tc_n}\leq 2+\displaystyle\min_{e\in\Ec}|\Fc_e|.
$$

If $\min\sigma(\Lc_1)=n$ alors $\lambda_{\Tc_n}\geq n$
and $n\leq 2+\displaystyle\min_{e\in\Ec}|\Fc_e|.$ Then the triangulation $\Tc_n$
is complete. By Proposition \ref{pro0}, $\lambda_{\Tc_n}=n.$
\end{dem}
\begin{rem}
In a homogeneous triangulation of a complete graph, we have always:
$$\lambda_{\Tc_n}=\min\sigma\left(\Lc^+_{1\scriptscriptstyle{\vert\ker(\delta^0)}}\right)=\min\sigma(\Lc_1).
$$
\end{rem}
\begin{figure}[!ht]
\centering
\begin{minipage}[t]{10cm}
\centering
\includegraphics*[width=13cm,height=9cm]{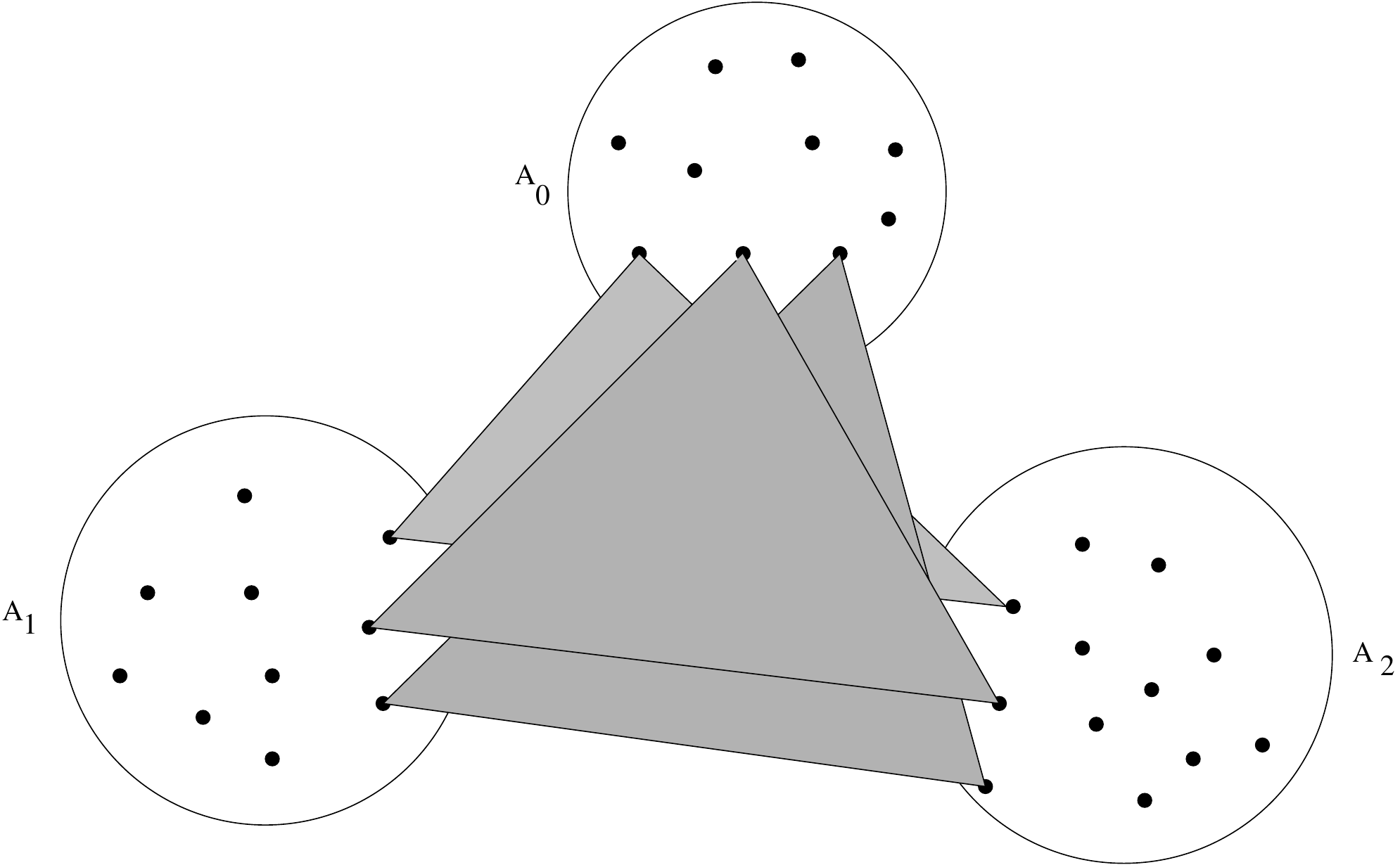}
\caption{Tripartite of a triangulation.}
\end{minipage}
\end{figure}

Our goal now is to find the best estimate of the eigenvalue $\lambda_{\Tc_n},$
when $\displaystyle\min_{e\in\Ec}|\Fc_e|$ is large enough.
Given $\Tc_n$ a homogeneous triangulation of a complete graph, \emph{the Cheeger constant}
$h(\Tc_n)$ is defined as follows, see (\cite{Chee}, \cite{PRT}):
$$h(\Tc_n):=\displaystyle\min_{\displaystyle\Vc=\displaystyle\bigcup_{i=0}^{2}\Ac_i}
\dfrac{n|\Fc(\Ac_0,\Ac_1,\Ac_2)|}{|\Ac_0||\Ac_1||\Ac_2|},
$$
where $\Ac_0,\Ac_1,\Ac_2$ are nonempty sets, considered as a partition of $\Vc,$ and $\Fc(\Ac_0,\Ac_1,\Ac_2)$
denotes the set of the oriented faces with one vertex in each $\Ac_i.$ It is clear that:
\begin{equation}\label{Hyp}
\exists x\in\Vc,\mbox{ such that }\forall e\in\Ec,~x\in\Fc_e\Longrightarrow h(\Tc_n)\neq0.
\end{equation}

Let $\Ac_0,\Ac_1\in\Vc$ are two finite sets such that $\Ac_0\cap\Ac_1=\emptyset,$ we define the set:
$$\Ec(\Ac_0,\Ac_1):=\left\{e\in\Ec,~\{e^-,e^+\}\cap\Ac_0\neq\emptyset
\mbox{ and }\{e^-,e^+\}\cap\Ac_1\neq\emptyset\right\}.
$$

\begin{thm}\label{thm}
If $h(\Tc_n)\neq0$ holds, then the spectral gap satisfies the following upper estimate
$$\lambda_{\Tc_n}\leq h(\Tc_n).
$$
\end{thm}
\begin{dem}
Set $\Ac_0,\Ac_1$ and $\Ac_2$ be a partition of $\Vc$ which realizes the minimum in $h(\Tc_n).$
We define $ \tilde{\psi}\equiv\tilde{\psi}(\Ac_0,\Ac_1,\Ac_2)\in l^2(\Ec)$ by
$$\tilde{\psi}(e)=\left\{\begin{array}{ll}
                        |\Ac_{i}| & \hbox{if }~ e\in\Ac_{\pi(i)}\times\Ac_{\pi(i+1)}\\
                        -|\Ac_{i}| & \hbox{if  }~e\in\Ac_{\pi(i+1)}\times\Ac_{\pi(i)}\\
                        0 & \hbox{else} \,\mbox{ \textit{i.e.} if }e\in\Ac_i,\mbox{ for all }i\in\{0,1,2\}

                 \end{array}
            \right.
$$
where $\pi:\Z\rightarrow\{0,1,2\}$ is a map such that $\pi(i)=i+1\mbox{ modulo }3.$
It is clear that $\tilde{\psi}\in\ker(\delta^0).$ Therefore, we have
\begin{equation*}
\begin{split}
\langle\Lc_1^+\tilde{\psi},\tilde{\psi}\rangle_{l^2(\Ec)}&=\|d^1\tilde{\psi}\|^2_{l^2(\Fc)}\\
&=\dfrac{1}{2}\displaystyle\sum_{(x,y,z)\in\Fc}
\left|\tilde{\psi}(x,y)+\tilde{\psi}(y,z)+\tilde{\psi}(z,x)\right|^2\\
&=n^2|\Fc(\Ac_0,\Ac_1,\Ac_2)|
\end{split}
\end{equation*}
and
\begin{equation*}
\begin{split}
\|\tilde{\psi}\|_{l^{2}(\Ec)}^{2}&=\dfrac{1}{2}\displaystyle\sum_{e\in\Ec}|\tilde{\psi}(e)|^2\\
&=\dfrac{1}{2}\left(|\Ec(\Ac_0,\Ac_1)||\Ac_2|^2+|\Ec(\Ac_1,\Ac_2)||\Ac_0|^2+|\Ec(\Ac_2,\Ac_0)||\Ac_1|^2\right)\\
&=|\Ac_0||\Ac_1||\Ac_2|^2+|\Ac_1||\Ac_2||\Ac_0|^2+|\Ac_2||\Ac_0||\Ac_1|^2\\
&=|\Ac_0||\Ac_1||\Ac_2|\left(|\Ac_0|+|\Ac_1|+|\Ac_2|\right)\\
&=n|\Ac_0||\Ac_1||\Ac_2|.
\end{split}
\end{equation*}

By the Rayleigh's principle, we have:
$$\lambda_{\Tc_n}=\min\sigma\left(\Lc^+_{1\scriptscriptstyle{\vert\ker(\delta^0)}}\right)
=\inf_{\psi\in\Cc_{c}(\Ec)\backslash\{0\}}\dfrac{\langle\Lc^+_1\psi,\psi\rangle_{l^{2}(\Ec)}}
{\|\psi\|_{l^{2}(\Ec)}^{2}}
\leq h(\Tc_n).
$$
\end{dem}
\begin{exa}
We consider the triangulation $\Tc_{4}=(\Kc_{4},\Fc)$ with $\Vc=\{1,2,3,4\}$
and $\Fc=\{(\textcolor{red}{1},2,3),(\textcolor{red}{1},2,4),(\textcolor{red}{1},3,4)\}.$
Then, we have $h(\Tc_{4})=2.$ Applying Theorem \ref{thm}, we obtain that
$$\lambda_{\Tc_4}\leq2.
$$
\end{exa}
\begin{exa}
We consider the triangulation $\Tc_{5}=(\Kc_{5},\Fc)$ with $\Vc=\{1,2,3,4,5\}$
and $\Fc=\{(\textcolor{red}{1},2,3),(\textcolor{red}{1},2,4),(\textcolor{red}{1},3,4),
(\textcolor{red}{1},2,5),(\textcolor{red}{1},3,5),(\textcolor{red}{1},4,5)\}.$
Then, we have $h(\Tc_{5})=\dfrac{5}{3}.$ Applying Theorem \ref{thm}, we obtain that
$$\lambda_{\Tc_5}\leq\dfrac{5}{3}.
$$
\end{exa}

\subsubsection{Lower estimates of spectral gap}
First, we recall some lower  estimates of the spectral gap of $\Lc_0$ on connected graphs,
see (\cite{CdV},\cite{KL}). Our objective is to find a lower estimate of the spectral gap
$\lambda_{\Tc_n}$ from the second eigenvalue of $\Lc_0.$
\begin{figure}[!ht]
\centering
\begin{minipage}[t]{10cm}
\centering
\includegraphics*[width=13cm,height=9cm]{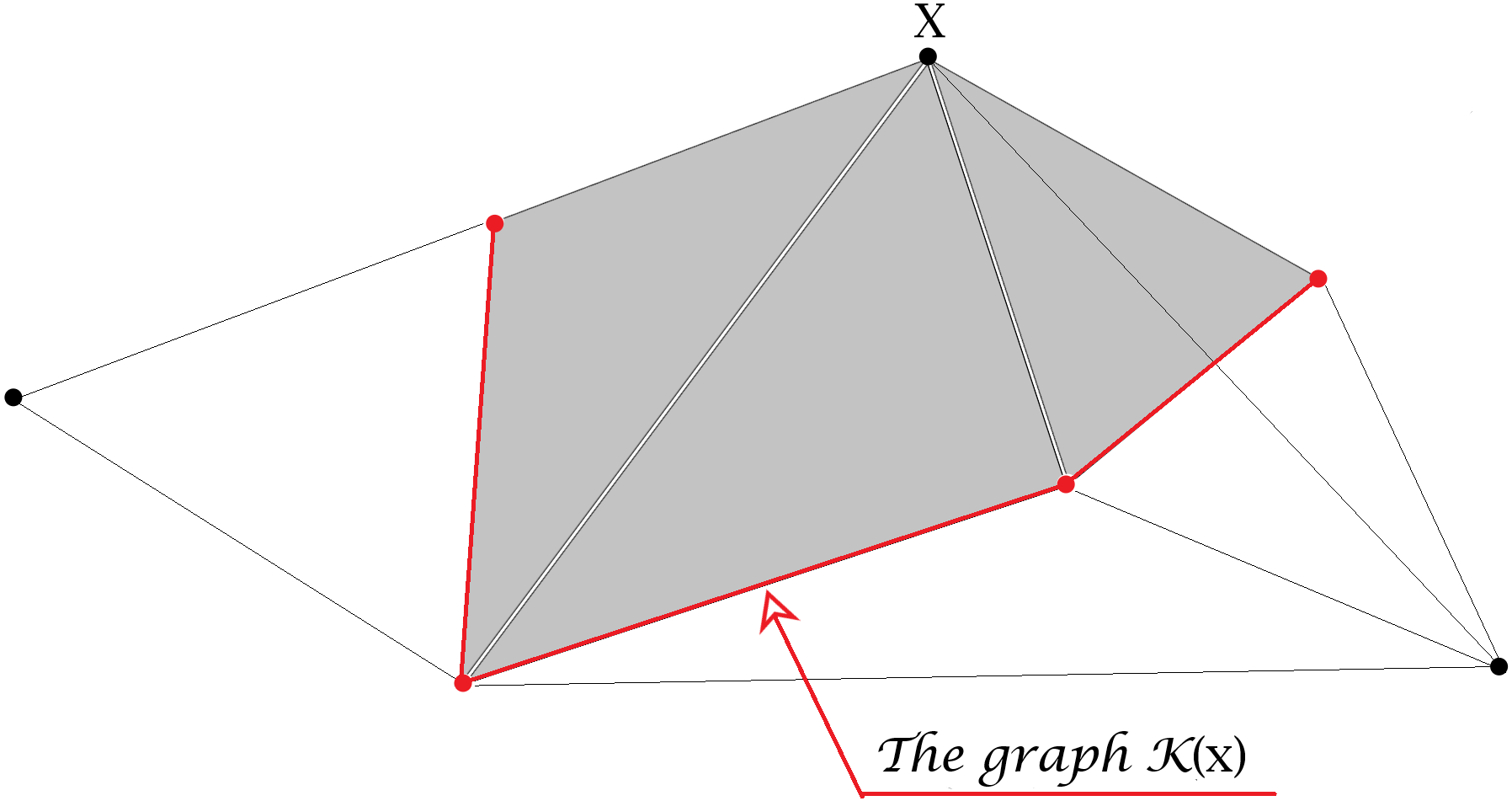}
\caption{The associated graph $\Kc$(x) of vertex x.}
\end{minipage}
\end{figure}

Let $\Tc=(\Vc,\Ec,\Fc)$ be a finite triangulation and $x\in\Vc.$
We define the set of edges $\Ec_{x}$:
$$\Ec_{x}:=\{e\in\Ec;~x\in\Fc_{e}\}.
$$

Denote by $\Kc(x):=(\Vc_{x},\Ec_{x})$ the associated graph of the vertex $x,$ where:
$$\Vc_{x}:=\{y\in\Vc;~\exists z\in\Vc\mbox{ such that }(y,z)\in\Ec_{x}\}.
$$

\begin{rem}
In a triangulation of a complete graph, the hypothesis \emph{(\ref{Hyp})}
assures that, for all $e\in\Ec,~\Fc_{e}\neq\emptyset.$ And then, we have that for all $x\in\Vc,~\Vc_x=\Vc\backslash\{x\}.$
\end{rem}

\begin{thm}\label{hamda}
Let $\Tc_n$ be a homogeneous triangulation of a complete graph.
Assume that, for all $e\in\Ec,~\Fc_{e}\neq\emptyset,$ then
$$\lambda_{\Tc_n}\geq\displaystyle\min_{e\in\Ec}\left(2\lambda_{1}\left(\Kc(e^-)\right)-|\Fc_{e}|\right).
$$
\end{thm}
\begin{dem}
Given $\varphi\in\ker(\delta^0)$ an eigenfunction associated of $\lambda_{\Tc_n}$ such that $\|\varphi\|_{l^2(\Ec)}=1.$ Then, we have
\begin{equation*}
 \begin{split}
    \lambda_{\Tc_n} & =\langle \Lc_{1}^{+}\varphi,\varphi\rangle_{l^2(\Ec)}\\
      & =\dfrac{1}{2}\displaystyle\sum_{e\in\Ec}\varphi(e)\displaystyle\sum_{x\in\Fc_{e}}
      (\varphi(e)+\varphi(e^+,x)+\varphi(x,e^-))\\
      & =\dfrac{1}{2}\displaystyle\sum_{e\in\Ec}|\Fc_{e}|\varphi^2(e)
      +\dfrac{1}{2}\displaystyle\sum_{e\in\Ec}\varphi(e)\displaystyle\sum_{x\in\Fc_{e}}
      (\varphi(e^+,x)+\varphi(x,e^-))\\
      & =\dfrac{1}{2}\displaystyle\sum_{e\in\Ec}|\Fc_{e}|\varphi^2(e)
      +\displaystyle\sum_{e\in\Ec}\varphi(e)\left(\displaystyle\sum_{x\in\Fc_{e}}
      \varphi(e^+,x)\right).
 \end{split}
\end{equation*}

For each vertex $x\in\Vc,$ we can define on $\Vc_{x}$ the function $f_{x}(y):=\varphi(x,y).$
Then, we have $f_{x}(y)=-f_{y}(x),$ for all $(x,y)\in\Ec.$ Next, we have
\begin{equation*}
 \begin{split}
 \Ic & =\displaystyle\sum_{x\in\Vc}\displaystyle\sum_{e\in\Ec_{x}}\varphi(e)\varphi(e^+,x)\\
 & =\displaystyle\sum_{x\in\Vc}\displaystyle\sum_{y\in\Vc_{x}}\displaystyle\sum_{e\in\Ec_{x},e^+=y}
 \varphi(e)\varphi(e^+,x)\\
 & =\displaystyle\sum_{x\in\Vc}\displaystyle\sum_{y\in\Vc_{x}}\displaystyle\sum_{e\in\Ec_{x},e^+=y}
 f_{e^-}(e^+)f_{e^+}(x)\\
 & =\displaystyle\sum_{x\in\Vc}\displaystyle\sum_{y\in\Vc_{x}}f_{x}(y)
 \left(\displaystyle\sum_{z\in\Fc_{(x,y)}}f_{y}(z)\right)\\
 & =\displaystyle\sum_{x\in\Vc}\displaystyle\sum_{y\in\Vc_{x}}f_{x}(y)
 \left(\displaystyle\sum_{z\in\Fc_{(x,y)}}f_{x}(y)-f_{x}(z)\right)\\
 & +\displaystyle\sum_{x\in\Vc}\displaystyle\sum_{y\in\Vc_{x}}f_{x}(y)
 \left(\displaystyle\sum_{z\in\Fc_{(x,y)}}f_{y}(x)+f_{y}(z)+f_{x}(z)\right).
 \end{split}
\end{equation*}

Thus,
$$\Ic =\displaystyle\sum_{x\in\Vc}\langle\Lc_{0\scriptscriptstyle{\Vc_{x}}}f_{x},f_{x}\rangle_{l^2(\Vc_{x})}
-2\langle deg_{\Ec}\varphi,\varphi\rangle_{l^2(\Ec)}
+\underbrace{\displaystyle\sum_{x\in\Vc}\displaystyle\sum_{y\in\Vc_{x}}f_{x}(y)
\left(\displaystyle\sum_{z\in\Fc_{(x,y)}}f_{y}(z)+f_{x}(z)\right)}_{\Jc}.
$$

Since $\Ec$ is the set oriented edges, we have
\begin{equation*}
 \begin{split}
 \displaystyle\sum_{x\in\Vc}\displaystyle\sum_{y\in\Vc_{x}}
 \displaystyle\sum_{z\in\Fc_{(x,y)}}\varphi(x,y)\varphi(y,z)
 & =\displaystyle\sum_{(x,y,z)\in\Fc}\varphi(x,y)\varphi(y,z)\\
 & =-\displaystyle\sum_{(x,y,z)\in\Fc}\varphi(x,y)\varphi(x,z)
 \end{split}
\end{equation*}

Then $\Jc=0.$ On other hand, we have
$$\varphi\in\ker(\delta^0)\Rightarrow\forall x\in\Vc,~f_{x}\perp1.
$$

Applying the Min-Max principle, we obtain that
\begin{equation*}
 \begin{split}
\lambda_{\Tc_n}&\geq\displaystyle\sum_{x\in\Vc}\lambda_{1}(\Kc(x))\|f_{x}\|^2_{l^2(\Vc_{x})}
-\langle deg_{\Ec}\varphi,\varphi\rangle_{l^2(\Ec)}\\
&=\displaystyle\sum_{x\in\Vc}\lambda_{1}(\Kc(x))\left(\displaystyle\sum_{y\in\Vc_{x}}f_{x}^2(y)\right)
-\dfrac{1}{2}\displaystyle\sum_{e\in\Ec}|\Fc_{e}|\varphi^2(e)\\
&=\dfrac{1}{2}\displaystyle\sum_{e\in\Ec}(2\lambda_{1}(\Kc(e^-))-|\Fc_{e}|)\varphi^2(e)\\
&\geq\displaystyle\min_{e\in\Ec}\left(2\lambda_{1}(\Kc(e^-))-|\Fc_{e}|\right).
 \end{split}
\end{equation*}
\end{dem}

\begin{rem}
This theorem is interesting in the case of triangulations of a complete graph
where its sub-graphs $(\Kc(x))_{x\in\Vc}$ are connected in such a way that
we have $\displaystyle\min_{x\in\Vc}\lambda_{1}(\Kc(x))\neq0.$
To ensure this assumption, we should take triangulations with a large number
of triangular faces.
\end{rem}

The next result is obtained from the lower estimate of the spectral gap
of the discrete Laplacian $\Lc_{0},$ see \cite{CGY} and \cite{CdV}.

\begin{coro}
Let $\Tc_n$ be a homogeneous triangulation of a complete graph.
Assume that for all $e\in\Ec,~\Fc_{e}\neq\emptyset,$ then
$$\lambda_{\Tc_n}\geq\displaystyle\min_{e\in\Ec}\left(\dfrac{h^2(\Kc(e^-))}{d_{e^-}}-|\Fc_{e}|\right).
$$
with for all $x\in\Vc,~d_{x}=\displaystyle\max_{y\in\Vc_{x}}deg(y).$
\end{coro}
\begin{dem}
Let $\varphi\in\ker(\delta^0)$ an eigenfunction of $\lambda_{\Tc_n}$ such that $\|\varphi\|_{l^2(\Ec}=1.$ Using Theorem \ref{hamda}, we obtain that
\begin{equation*}
  \begin{split}
     \lambda_{\Tc_n} & \geq\displaystyle\min_{e\in\Ec}\left(2\lambda_{1}(\Kc(e^-))-|\Fc_{e}|\right)\\
                     & \geq\displaystyle\min_{e\in\Ec}\left(\dfrac{h^2(\Kc(e^-))}{d_{e^-}}-|\Fc_{e}|\right).
  \end{split}
\end{equation*}
\end{dem}

\begin{exa}
We consider the triangulation $\Tc_5=(\Kc_5,\Fc)$ with $\Vc=\{1,2,3,4,5\}$
and $\Fc=\{(1,2,3),(1,2,5),(1,3,4),(1,4,5),(2,3,4),(2,4,5)\}.$
Using Theorem \ref{hamda}, we obtain that:
$$\lambda_{\Tc_5}\geq2\lambda_{1}(\Cc_{4})-2=2\times2-2=2.
$$
\end{exa}

\begin{exa}
We consider the triangulation $\Tc_6=(\Kc_6,\Fc)$ with $\Vc=\{1,2,3,4,5,6\}$ and
$\Fc=\{(1,2,3),(1,2,6),(1,3,4),(1,4,5),(1,5,6),(2,3,5),(2,4,5),(2,4,6),(3,4,6),(3,5,6)\}.$
Using Theorem \ref{hamda}, we obtain that: $$\lambda_{\Tc_6}\geq2\lambda_{1}(\Cc_{5})-2\approx2\times1.3819-2=0.7639.
$$
\end{exa}

\begin{center}
\begin{figure}[!ht]
\begin{minipage}[t]{10cm}
\centering
\includegraphics*[width=13cm,height=7cm]{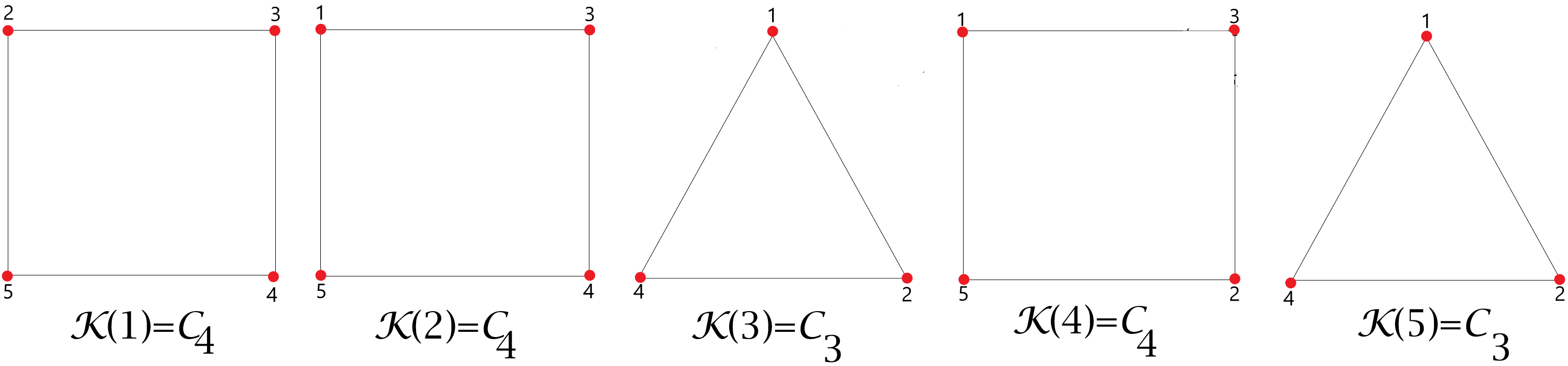}
\caption{The sub-graphs of the triangulation $\Tc_{5}$}
\end{minipage}
\end{figure}
\end{center}

\begin{figure}[!ht]
\begin{minipage}[t]{10cm}
\centering
\includegraphics*[width=13cm,height=7cm]{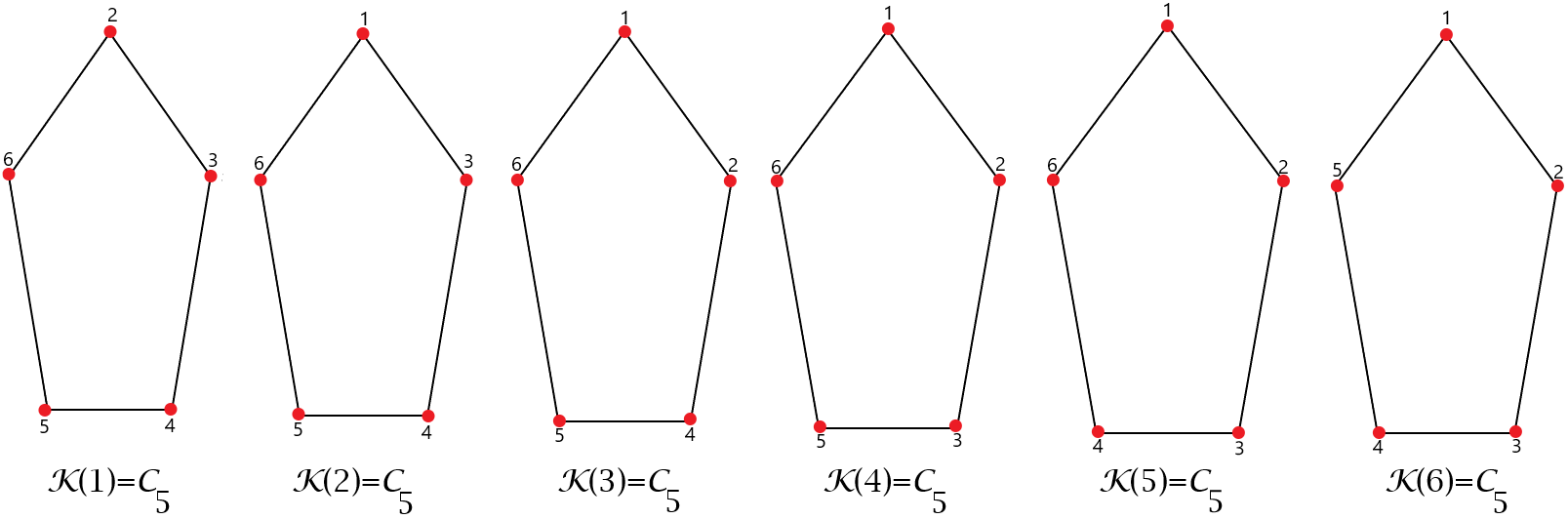}
\caption{The sub-graphs of the triangulation $\Tc_{6}$}
\end{minipage}
\end{figure}

The spectral gap of a discrete graph is a monotonously increasing function
of the set of edges. In other words, adding an edge always increases of
the second eigenvalue or keeps it unchanged, provided that we have the same set of
vertices, see \cite{KMN}.
\begin{propr}\label{exxp}
Let $\Kc$ be a connected graph and let $\Kc^{'}$ be a graph obtained from $\Kc$
by adding one edge between two vertices. Then the following hold:
$$\lambda_{1}(\Kc)\leq\lambda_{1}(\Kc^{'}).
$$
\end{propr}

\begin{exa}
We consider the triangulation $\Tc^{'}_5=(\Kc_5,\Fc)$ with $\Vc=\{1,2,3,4,5\}$
and $\Fc=\{(1,2,3),(1,2,5),(1,3,4),(1,3,5),(1,4,5),(2,3,4),(2,4,5),(3,4,5)\}.$
Using Proposition \ref{exxp}, we obtain that $\displaystyle\min_{i=1,...,5}\lambda_{1}(\Kc(i))=\lambda_{1}(\Cc_{4})=2.$
By Theorem \ref{hamda}, we have:
$$\lambda_{\Tc^{'}_5}\geq2\lambda_{1}(\Cc_{4})-3=2\times2-3=1.
$$

\begin{figure}[!ht]
\begin{minipage}[t]{10cm}
\centering
\includegraphics*[width=13cm,height=8cm]{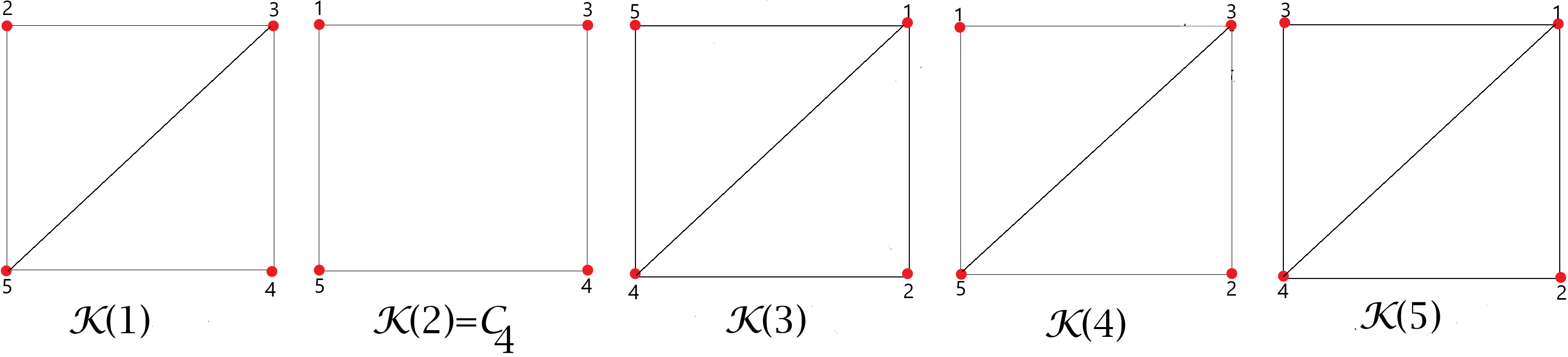}
\caption{The sub-graphs of the triangulation $\Tc^{'}_{5}$}
\end{minipage}
\end{figure}

\textbf{\textit{\emph{Acknowledgments}}:} \emph{I would like to sincerely thank my PhD advisors, Professors Colette Ann\'e and Nabila Torki-Hamza for helpful discussions. I am very thankful
to them for all the encouragement, advice and inspiration. I would also like to thank
the Laboratory of Mathematics Jean Leray and the research unity (UR/13 ES 47) for their
continuous support. This work was partially financially supported by the "PHC Utique"
program of the French Ministry of Foreign Affairs and Ministry of higher education
and research and the Tunisian Ministry of higher education and scientific research
in the CMCU project number 13G1501 "Graphes, G\'{e}om\'{e}trie et Th\'{e}orie Spectrale".}


\begin{thebibliography}{xxxxxx}
\bibitem[1]{Ay} H. Ayadi: \emph{Spectra of Laplacians on an infinite graph},
Oper. Matrices, 11, no.2, 567-586, 2017.
\bibitem[2]{AM} N. Alon and V.D. Milman: $\lambda_{1},$ \emph{isoperimetric inequalities for graphs,
and superconcentrators}, J. Combin. Theory Ser. B, 38, no. 1, 73-88, 1985.
\bibitem[3]{AT} C. Ann\'{e} and N. Torki-Hamza: \emph{The Gau\ss -Bonnet operator
of an infinite graph}, Anal. Math. Phys. 5, no.2, 137-159, 2015.
\bibitem[4]{Che} Y. Chebbi: \emph{The discrete Laplacian of a 2-simplicial complex},
Potential Analysis, 49, no.2, 331-358, 2018.
\bibitem[5]{Chee} Y. Chebbi: \emph{Laplacien discret d'un 2-complexe simplicial},
HAL Id: tel-01800569, 2018.
\bibitem[6]{Ch} J. Cheeger: \emph{A lower bound for the smallest eigenvalue of the Laplacian},
Problems in Analysis (Papers dedicated to Salamon Bochner, 1969), pp. 195-199,
Princeton University Press, Princeton, NJ, 1970.
\bibitem[7]{Chu} F.R.K Chung: \emph{Spectral graph theory}, CBMS Regional
Conference Series in Mathematics. 92. Providence, RI: American Mathematical
Society (AMS). xi, p. 207, 1994.
\bibitem[8]{CGY} F. Chung, A. Grigoryan, S-T. Yau: \emph{Higher eigenvalues and isoperimetric
inequalities on Riemannian manifolds and graphs}, Comm. Anal. Geom., 8, no. 5, 969-1026, 2000.
\bibitem[9]{CdV} Y. Colin de Verdi\`{e}re: \emph{Spectres de grahes}, Cours Sp\'{e}cialis\'{e}s,
4. Soci\'{e}t\'{e} Math\'{e}matique de France, Paris, 1998.
\bibitem[10]{D} J. Dodziuk: \emph{Difference equations, isoperimetric inequality and transience
of certain random walks}, Trans. Amer. Math. Soc. 284, 787-794, 1984.
\bibitem[11]{Eck} B. Eckmann, \emph{Harmonische funktionen und randwertaufgaben in einemkomlex},
Commentarii Mathematici Helvetici, vol 17, no. 1, 240-255, 1944.
\bibitem[12]{FP} D. L. Ferrario et R. A. Piccinini:
\emph{Simplicial Structures in Topology}, CMS Books in Mathematics, p. 243, 2011.
\bibitem[13]{G} S. Gol\'{e}nia: \emph{Hardy inequality and asymptotic eigenvalue distribution
for discrete Lapacians}, Journal of Functional Analysis, vol. 266, no.5, 2662-2688, 2014.
\bibitem[14]{KL} M. Keller and D. Lenz: \emph{Unbounded laplacians on graphs:basic spectral
properties and the heat equation}, Math. Model. Nat. Phenom.5, no. 4, 198-224, 2010.
\bibitem[15]{KMN} P. Kurasov, G. Malenov\'{a} and S. Naboko:
\emph{Spectral gap for quantum graphs and their connectivity}, Journal of Physics A:
Mathematical and Theoretical. vol 46, no. 27, 2013.
\bibitem[16]{PRT} O. Parzanchevski, R. Rosenthal and R. J. Tessler: \emph{Isoperimetric inequalities
in simplicial complexes}, COMBINATORICA, vol. 36, no. 2, 195-227, 2016.
\bibitem[17]{RS} M. Reed , B. Simon: \emph{Methods of Modern Mathematical Physics},
Vol. 1, New York Academic Press, 1980.
\end{thebibliography}
\end{document}